\documentclass[draft]{birkmult}

\usepackage{amsmath, amsthm, amssymb, mathrsfs}

\newcommand{\D}    {\mathbb{D}}

\newcommand{\C}    {\mathbb{C}}
\newcommand{\N}    {\mathbb{N}}

\newcommand{\pl}   {\mathcal{P}}

\newcommand{\Ml}   {\mathcal{M}}
\newcommand{\Nl}   {\mathcal{N}}

\newcommand{\supp}  {\textnormal{supp}}

\newcommand{\cp}      {\textnormal{cap}}

\newcommand{\cws}  {\stackrel{*}{\to}}

\numberwithin{equation}{section}
 
\newtheorem{theorem}{Theorem}[section]
\newtheorem{lemma}[theorem]{Lemma}

\newtheorem{example}{Example}[section]
 \theoremstyle{remark}
\newtheorem*{rem}{Remarks}

\begin{document}

\title[Ratios of Norms for Polynomials and Connected $n$-width Problems]{Ratios of Norms for Polynomials and\\ Connected $n$-width Problems}
%----------Author 1
\author{V. A. Prokhorov}

\address{Department of Mathematics and Statistics \\
University of South Alabama \\
Mobile, Alabama 36688-0002, USA}

\email{prokhoro@jaguar1.usouthal.edu}

%----------Author 2
\author{E. B. Saff}
\address{Department of Mathematics \\
Vanderbilt University \\
Nashville, TN 37240, USA}

\email{edward.b.saff@vanderbilt.edu}

%----------Author 3
\author{M. Yattselev}
\address{INRIA, Project APICS\newline
2004 route des Lucioles --- BP 93\newline
06902 Sophia-Antipolis, France}

\email{myattsel@sophia.inria.fr}

\thanks{The research of E. B. Saff was supported, in part, by U.S. National Science Foundation under grant DMS-0603828. The research of M. Yattselev was supported, in part, by U.S. National Science Foundation under grant DMS-0603828 (E. B. Saff principal investigator).}



\date{April 23, 2008}
\dedicatory{Dedicated to Bjorn Gustafsson on the occasion of his 60th birthday}
%%% ----------------------------------------------------------------------

\begin{abstract}
Let $G\subset  \C$ be a bounded simply connected domain  with boundary $\Gamma$  and let $E\subset G$ be a regular compact set with  connected complement. In this paper we investigate  asymptotics of the extremal constants:
\[
\chi_n = \inf_{p\in\pl_{k_n}} \sup_{q\in\pl_{n-k_n}}\frac{||pq||_E}{||pq||_\Gamma},\quad n=1, 2, \dots,
\]
where $\|\cdot\|_K$ is the supremum norm on a compact set $K$, $\pl_m$ is the set of all algebraic polynomials of degree at most $m$, and $k_n/n\to\theta\in[0,1]$ as $n\to\infty$. Subsequently, we obtain asymptotic behavior of  the Kolmogorov $k$-widths, $k=k_n$,  of  the unit ball $A_n^\infty$ of $H^\infty\cap\pl_n$ restricted to $E$ in $C(E)$, where $H^\infty$ is the Hardy space of bounded analytic functions on $G$ and  $C(E)$ is the space of continuous functions on $E$.
\end{abstract}

\maketitle

\section{Introduction}
\label{sec:intro}

The Kolmogorov $k$-width of a set $A$ contained in a Banach space $X$ is defined by
\[
d_k(A;X) := \inf_{X_k}\sup_{h\in A}\inf_{g\in X_k}\|h-g\|,
\]
where $X_k$ runs over all $k$-dimensional subspaces of $X$ and $\|\cdot\|$ is a norm on $X$. Let $G$ be a bounded simply connected domain with boundary $\Gamma$ in the complex plane $\C$, and $H^\infty$ be the Hardy space of bounded analytic functions in $G$. Denote by  $E\subset G$  a regular compact set with connected complement $D$ and $A^\infty$  the unit ball of $H^\infty$ restricted to $E$. In \cite{Wid72} H. Widom investigated the asymptotic behavior of $d_k(A^\infty;C(E))$, where $C(E)$ is the space of continuous functions on $E$ endowed with the usual supremum norm $\|\cdot\|_E$. It is proved that
\begin{equation}
\label{eq:kWidthLim}
\lim_{k\to\infty}\left(\frac1k\log d_k(A^\infty;C(E))\right) = -\frac{1}{\cp(E,\Gamma)},
\end{equation}
where $\cp(E,\Gamma)$ is the {\it condenser (Green) capacity} of $E$ with respect to $G$ (see, for example,  \cite[Sec. II.5]{SaffTotik}). Further, in \cite{FishMic80} (see also \cite[Sec. 7.5]{Fisher}) S. D. Fisher and C. A. Micchelli obtained the following representation for $d_k(A^\infty;C(E)):$
\begin{equation}
\label{eq:kWidthHol}
d_k(A^\infty;C(E)) = \inf_{z_1,\ldots,z_k}\sup\{\|h\|_E:~h\in A^\infty, \;\; h(z_j)=0, \;\; j=1,\ldots,k\}.
\end{equation}
Clearly, it is enough to consider only the {\it Blaschke products} instead of all functions from $A^\infty$ in (\ref{eq:kWidthHol}). Then it is a consequence of \cite{FishS99} that the zero counting measures of any asymptotically extremal sequence of Blaschke products swept out to $\partial E$ converge weak-star to the {\it Green equilibrium distribution} on $E$ relative to $G$.

In this paper we investigate the $n$-th root behavior of $d_k(A_n^\infty;C(E))$, $k=k_n$, the Kolmogorov $k$-widths of the unit ball $A_n^\infty$ of $H^\infty\cap\pl_n$ restricted to $E$ in $C(E)$, and show its connection to the following extremal problem:
\begin{equation}
\label{eq:ExPr1}
\chi_n = \inf_{p\in\pl_{k_n}} \sup_{q\in\pl_{n-k_n}} \frac{||pq||_E}{||pq||_\Gamma},
\end{equation}
where
\begin{equation}
\label{eq:theta}
\lim_{n \to \infty} \frac{k_n}{n}=\theta, \;\;\; \theta\in[0,1].
\end{equation}
Analogous $k$-width occur in the study of truncated Hankel operators which the authors will explore in a later paper.

Regarding the minimax problem defined in (\ref{eq:ExPr1}), we observe that it connects two well-understood extremal problems of potential theory. It is an simple consequence of the Bernstein-Walsh inequality (\cite{Walsh} and \cite[Sec. III.2]{SaffTotik}) and properties of the Chebyshev polynomials for $E$ that
\begin{equation}
\label{eq:ExPrOnE}
\lim_{n\to\infty}\left(\inf_{p\in\pl_n}\frac{\|p\|_E}{\|p\|_\Gamma}\right)^{1/n} = \exp\left\{-\max_{z\in\Gamma}g(z,\infty)\right\},
\end{equation}
where $g(\cdot,\infty)$ is the Green function for $D$ with singularity  at infinity. It is easy to see that
\begin{equation}
\label{eq:ExPrOnGamma}
\lim_{n\to\infty}\left(\sup_{q\in\pl_n}\frac{\|q\|_E}{\|q\|_\Gamma}\right)^{1/n} = 1
\end{equation}
 and  the extremal polynomial is  $q\equiv1$. Furthermore, it is readily verified that polynomials $z^n-R^n$ are asymptotically extremal for (\ref{eq:ExPrOnGamma}) whenever $R$ is such that  $\{|z|<R\}\supset G$. Let us also illustrate extremal problem (\ref{eq:ExPrOnE}). Put $E$ to be the closed unit disk $\overline\D$ and $\Gamma$ to be the circle of radius $R>2$ centered at 1. In this case $g(z,\infty)=\log|z|$ and therefore the monomials $z^n$ are extremal for (\ref{eq:ExPrOnE}) and the limit is equal to $1/(R+1)$. Moreover, the polynomials $z^n-1$ are asymptotically extremal for that problem. 

This paper is organized as follows. In Section \ref{sec:pot} we consider two minimal energy problems, one for the Green potentials and another for the logarithmic potentials, that are vital for our main results. The latter are given in  Section \ref{sec:main}, which  contains results on the behavior of $\chi_n$ and the extremal polynomials (Theorems \ref{thm:1} and \ref{thm:2}) as well as connection with $n$-width (Theorem \ref{thm:3}). In Section \ref{sec:ContExtrPr} we study some extremal problem of the potential theory which can be considered as an continuous analog of the extremal problem (\ref{eq:ExPr1}). In Section \ref{sec:lm} we provide a detailed description of the extremal measures defined in Section \ref{sec:pot}. Sections \ref{sec:proofs1} and \ref{sec:proofs2}  of this paper consist of proofs of the stated results. In Section \ref{sec:proofs3} we investigate the asymptotics of $k$-widths.

\section {Equilibrium Measures}
\label{sec:pot}

Let $G$, $\Gamma$, $E$, and $D$ be as described. We shall use the standard terminology that a property holds {\it quasi-everywhere} (q.e.) if it holds everywhere except for a set of zero {\it logarithmic capacity} (see \cite[Sec. I.1]{SaffTotik} or \cite[Sec. 5]{Ransford} for the definition of capacity). In this paper we extensively utilize {\it logarithmic} and {\it Green potentials}. The logarithmic potential of a finite positive Borel measure $\nu$ with compact support $\supp(\nu)$, is given by
\[
U^\nu(z)=-\int\log|z-t|d\nu(t).
\]
It is superharmonic in $\C$ and  harmonic in $\C\setminus\supp(\nu)$. Unlike the logarithmic case, Green potentials are defined relative to a domain. Let $\nu$ be a  positive Borel measure compactly supported in $D$. Then the Green potential of $\nu$ relative to $D$ is given by
\[
U_D^\nu(z) =\int g(z,t)d\nu(t),
\]
where $g(z,t)=g_E(z, t)$ is the Green function for $D$ with singularity at $t\in D$. Since $E$ is a regular compact set, $g(z,t)=0$ for $z\in\partial D = \partial E$. Here and in what follows we assume that  $g(z,t)=0$  for all $z\in E$. The Green potential of $\nu$ is nonnegative and superharmonic in $D$,  harmonic in $D\setminus\supp(\nu)$, and satisfies $U_D^\nu=0$ on $E$. 
 
Let $K$ be a compact set.  Denote by $\Lambda_\delta(K)$, $\delta>0$, the set of positive Borel measures $\lambda$ of mass $\delta=|\lambda|= \int d \lambda$ compactly supported on $K$. 

For each $\theta\in[0,1)$ consider the following weighted Green energy of a measure $\lambda\in\Lambda_{1-\theta}(\Gamma)$:
\begin{equation}
\label{eq:JEnergy}
J_{\theta}(\lambda) := \int \int g(z, t) d \lambda (t) d \lambda(z) -2 \int g(t, \infty) d \lambda(t).
\end{equation}
Then we have the following result.
\begin{theorem}
\label{thm:aux1}
For each $\theta\in[0,1)$ there exists a unique measure $\lambda_\theta\in\Lambda_{1-\theta}(\Gamma)$  such that
\begin{equation}
\label{eq:m11}
 \min_{\lambda\in\Lambda_{1-\theta}(\Gamma)} J_{\theta}(\lambda) =J_{\theta}(\lambda_\theta).
\end{equation}
The extremal measure $\lambda_\theta$ satisfies the following properties:
\begin{equation}
\label{eq:m12}
U_D^{\lambda_\theta}(z) -g(z,\infty) = m_\theta, \;\;\; z\in S_\theta:=\supp(\lambda_\theta)\subseteq \Gamma,
\end{equation}
and
\begin{equation}
\label{eq:m10}
U_D^{\lambda_\theta}(z) -g(z,\infty) \geq m_\theta, \;\;\; z\in\Gamma,
\end{equation}
where
\begin{equation}
\label{eq:m15}
m_\theta := \frac{1}{1-\theta}\left(J_{\theta}(\lambda_\theta)+\int g(t,\infty)d\lambda_\theta(t)\right).
\end{equation}
\end{theorem}
\begin{rem}
\begin{itemize}
\item[(a)] This theorem is a special case of \cite[Thm. II.5.10]{SaffTotik} for the external field $-g(\cdot,\infty)/(1-\theta)$. We exhibit the dependence of $\lambda_\theta$ on $\theta$ in Theorem \ref{thm:4} (see also \cite[Thm. 2.4]{LevLub01}). In particular, $\lambda_{\theta_1}-\lambda_{\theta_2}$ is a positive measure for any choice of $\theta_1<\theta_2$.
\item[(b)] In general, (\ref{eq:m12}) holds only q.e. on $S_\theta$. However, as pointed out in \cite[Thm. 2.2]{LevLub01}, the regularity of $\Gamma$ is sufficient for this property to hold at every point of $S_\theta$.
\item[(c)] As shown later in Lemma \ref{lem:aux2}, $m_0=0$, $\lambda_0=\omega_\Gamma$, and $S_0=\Gamma$, where $\omega_K$ stands for the {\it logarithmic  equilibrium distribution} on a set $K$.
\item[(d)] It follows from Theorem \ref{thm:4}, the limit of $m_\theta$ as $\theta$ approaches $1$ from the left exists and
\begin{equation}
\label{eq:mone}
m_1:=\lim_{\theta\to1^-}m_\theta=-\max_{z\in\Gamma}g(z,\infty).
\end{equation}
Furthermore,   we define $\lambda_1$ to be the zero measure.
\item[(e)] The measure $\lambda_\theta$ is uniquely determined by conditions (\ref{eq:m12}) and (\ref{eq:m10}). If 
$\lambda \in \Lambda_{1-\theta}(\Gamma)$ has a finite Green energy, $U_D^\lambda-g(z, \infty)=c$ on $\supp (\lambda)$  and $U_D^\lambda -g(z, \infty)\ge c$ on $\Gamma$, then $\lambda =\lambda_\theta$ and $c=m_\theta$ (see \cite[Thm. II. 5.12]{SaffTotik}).
\end{itemize}
\end{rem}

Let us consider the special case when $\Gamma$  is a level curve of $g(\cdot,\infty)$. 
\begin{example}
\label{cor:1}
Let $\Gamma =  \{z:~g(z,\infty)=R\}$ for some constant $R>0$. Then for every $\theta\in[0,1)$ we have
\[
\lambda_\theta =  (1-\theta)\omega_\Gamma \;\;\; \mbox{and} \;\;\; m_\theta=-\theta R.
\]
\end{example}
It is easy to see (cf. Lemma \ref{lem:aux2}) that 
\[
U_D^{(1-\theta)\omega_\Gamma}(z)-g(z,\infty)=-\theta g(z,\infty)=-\theta R, \;\;\; z\in\Gamma.
\]
Therefore, by Remark (e), we get
\[
\lambda_\theta = (1-\theta)\omega_\Gamma, \;\;\; m_\theta=-\theta R.
\]

The second  extremal problem that we need and which, in a way, is complimentary to (\ref{eq:m11}), is related to the following energy integral:
\[
I_\theta(\mu):=\int U^\mu(t) d\mu(t)+2\int U^{\lambda_\theta}(t)d\mu(t),
\]
where $\theta\in(0,1]$ and $\mu\in\Lambda_\theta(E)$.

As in the case of Theorem \ref{thm:aux1}, the following  is known \cite[Thm. I.1.3]{SaffTotik}.
\begin{theorem}
\label{thm:aux2}
For each $\theta\in(0,1]$ there exists a unique measure $\mu_\theta\in\Lambda_\theta(E)$ such that
\begin{equation}
\label{eq:m21}
\min_{\mu\in\Lambda_\theta(E)}I_\theta(\mu)=I_\theta(\mu_\theta).
\end{equation}
Moreover, the extremal measure $\mu_\theta$ has the following properties:
\begin{equation}
\label{eq:m22}
U^{\mu_\theta+\lambda_\theta}(z)=\widehat m_\theta, \;\;\; z\in\supp(\mu_\theta)\subseteq\partial E,
\end{equation}
and
\begin{equation}
\label{eq:m20}
U^{\mu_\theta+\lambda_\theta}(z) \geq \widehat m_\theta, \;\;\; z\in E,
\end{equation}
where
\begin{equation}
\label{eq:m25}
\widehat m_\theta := \frac1\theta\left(I_{\theta}(\mu_\theta)-\int U^{\lambda_\theta}(t)d\mu_\theta(t)\right).
\end{equation}
\end{theorem}
\begin{rem} 
\begin{itemize}
\item[(a)] For $\theta=1$, (\ref{eq:m21}) reduces to the classical (unweighted) minimal energy problem (cf. \cite[Sec. 3.3]{Ransford} and \cite[Sec. I.1]{SaffTotik}) when $\theta=1$. In this case, $\mu_1$ is the logarithmic  equilibrium distribution $\omega_E$ and and $\widehat m_1$ is the Robin constant for $E$, $\widehat m_1=-\log\cp(E)$, where $\cp(E)$ is the logarithmic capacity of $E$.
\item[(b)] It is a well-known fact that $\supp(\mu_\theta)\subseteq\partial E$ (see, for example,   \cite[Thm. IV.1.10(a)]{SaffTotik}).
\item[(c)] As in the case of $\lambda_1$, it is convenient for us to define $\mu_0$ to be the zero measure.
\item[(d)] The measure $\lambda_\theta$ is uniquely determined by conditions (\ref{eq:m22}) and (\ref{eq:m20}). If 
$\mu \in \Lambda_{\theta}(E)$ has a finite  energy and $U^{\mu+\lambda_\theta}=c$ on $\supp (\mu)$  and $U^{\mu+\lambda_\theta} \ge c$ on $E$, then $\mu=\mu_\theta$ and $c=\widehat{m_\theta}$ (cf.  \cite[Thm. I. 3.3]{SaffTotik}).
\end{itemize}
\end{rem}

Further properties of $\lambda_\theta$ and $\mu_\theta$ and the constants $m_\theta$ and $\widehat m_\theta$ are given in Section \ref{sec:lm}, including asymptotics as $\theta\to0$.

\section{Main Results}
\label{sec:main}

Let $\chi_n$ be defined by (\ref{eq:ExPr1}) and  (\ref{eq:theta}). Below  we show that $\lim \chi_n^{1/n}$ exists and provide the asymptotic behavior of the zeros of the {\it extremal polynomials}. The latter are defined as follows. Let $\{p_n, q_n\}_{n\in\N}$, $p_n\in\pl_{k_n}$ and $q_n\in\pl_{n-k_n}$, be such that
\begin{equation}
\label{eq:ExPoly}
\lim_{n\to\infty} \left(\frac{1}{\chi_n}\frac{\|p_nq_n\|_E}{\|p_nq_n\|_\Gamma}\right)^{1/n} = \lim_{n\to\infty}\left(\frac{1}{\chi_n}\sup_{q\in\pl_{n-k_n}}\frac{\|p_nq\|_E}{\|p_nq\|_\Gamma}\right)^{1/n} = 1.
\end{equation}
We shall call $\{p_n,q_n\}$, satisfying the equalities above, a {\it sequence of asymptotically extremal pairs} of polynomials. To each such $p_n$ and $q_n$ we associate the zero counting measures, $\nu(p_n)$ and $\nu(q_n)$, respectively, defined by the rule
\[
\nu(p_n):=\frac1n\sum_{p_n(z)=0}\delta_z \;\;\; \mbox{and} \;\;\; \nu(q_n):=\frac1n\sum_{q_n(z)=0}\delta_z,
\]
where $\delta_z$ is the point mass distribution at $z\in\C$ and the sums are taken counting multiplicities of zeros of polynomials 
$p_n$ and $q_n$.
It also will be convenient for us to sweep out ({\it balayage}) measures $\nu(p_n)$ and $\nu(q_n)$ onto $\partial E$ and $\Gamma$, respectively. Recall that for any finite positive Borel measure $\nu$ compactly supported in $\C$ and, with finite energy if $\supp(\nu)\cap\partial D\neq\emptyset$,   there exists a unique measure $\widehat\nu$, the balayage measure of $\nu$, supported on $\partial E$, such that $|\nu|=|\widehat\nu|$, 
\[
U^{\widehat\nu}(z) = U^\nu(z) + \int g(t,\infty)d\nu(t), \;\;\; z\in E, \;\;\; \mbox{if} \;\;\; \supp(\nu)\subset D,
\]
and
\[
U^{\widehat\nu}(z) = U^\nu(z), \;\;\; z\in\overline D, \;\;\; \mbox{if} \;\;\; \supp(\nu)\subseteq E.
\]
We remark that for any positive compactly supported in $D$  measure $\nu$,
\begin{equation}
U_D^\nu(z)=U^{\nu-\widetilde \nu}(z)+\int g(t, \infty) d \nu(t), \quad z \in \C.
\label{10.29.07.1}
\end{equation}
Denote by $\widetilde \nu$ the balayage of a finite positive Borel measure $\nu$ compactly supported in $\C\setminus\overline G$ onto $\Gamma$.
We have $|\widetilde \nu|=|\nu|$ and
\begin{equation}
U^{\widetilde \nu }(z)=U^\nu(z)+\int g_{\overline G }(t, \infty) \ d \nu(t), \quad z \in \overline G, 
\label{10.29.07.2}
\end{equation}
where $g_{\overline G}(z, \infty)$ is the  Green function of the domain $\overline \C\setminus \overline G$ with singularity at infinity.
 Now we define  measures $\alpha(p_n)$ and $\beta(q_n)$ as
\begin{equation}
\alpha(p_n) :=  \nu(p_n)_{|\overline D} + \widehat{\nu(p_n)_{|E\setminus\partial E}},
\label{11.03.07.1}
\end{equation}
and
\begin{equation}
\beta(q_n) := \nu(q_n)_{|\overline G} + \widetilde{\nu(q_n)_{|\C\setminus\overline G}} + \frac{n-k_n-\deg(q_n)}{n} \ \omega_\Gamma,
\label{11.03.07.2}
\end{equation}
respectively, where  a notation $\lambda _{| K}$ means restriction of a measure $\lambda$ on a set $K$.
Let $\Ml(E)=\{\nu:\nu \in \Lambda_\theta(E), \ \widehat \nu=\mu_\theta\}.$

 The following result holds.

\begin{theorem}
\label{thm:1} 
Let $\{k_n\}_{n\in\N}$ satisfy (\ref{eq:theta}) for some $\theta\in[0,1]$. Then
\begin{equation}
\label{eq:L11}
\lim_{n\to\infty}\left(\frac1n\log \chi_n\right) = m_\theta,
\end{equation}
where $m_\theta$ was defined in (\ref{eq:m15}) and (\ref{eq:mone}). If $\{p_n,q_n\}$ is a sequence of asymptotically extremal pairs of polynomials in the sense (\ref{eq:ExPoly}), then, for $\theta\in(0,1)$ any weak-star limit point of $\{\nu(p_n)\}$ belongs to $\Ml(E)$ and 
\begin{equation}
\label{eq:L12}
\alpha(p_n) \cws \mu_\theta \;\;\;  \mbox{as} \quad n \to \infty,
\end{equation}
where $\cws$ stands for the convergence of measures in the weak-star sense.
Moreover,
\begin{equation}
\label{eq:L13}
\begin{array}{ll}
\nu(q_n) \cws \lambda_\theta & \mbox{as} \  n \to \infty, \ \mbox{if} \;\;\; \overline\C\setminus S_\theta \;\;\; \mbox{is connected}, \\
\beta(q_n) \cws \lambda_\theta 
& \mbox{as } \  n \ \to \infty ,\ \mbox{otherwise}.
\end{array}
\end{equation}
\end{theorem}
\begin{rem}
\begin{itemize}
\item [(a)] Since $m_0=0$ and $m_1=-\max_{z\in\Gamma}g(z,\infty)$, (\ref{eq:ExPr1}) indeed connects extremal problems (\ref{eq:ExPrOnE}) and (\ref{eq:ExPrOnGamma}).
\item [(b)] Observe that $\alpha(p_n)=\nu(p_n)$ when $E$ has empty interior. In this case (\ref{eq:L12}) is a statement on the convergence of counting measures themselves, rather than their balayages onto $\partial E$.
\end{itemize}
\end{rem}

The following theorem is related to  the case when $k_n \to \infty$ and $k_n=o(n)$ as $n\to\infty$. To formulate the result, we need to slightly modify the definition of an asymptotically extremal sequence. We say that a sequence $\{p_n\}$ is {\it asymptotically extremal} if
\begin{equation}
\label{eq:ExPoly1}
\lim_{n\to\infty}\left(\frac{1}{\chi_n} \sup_{q\in\pl_{n-k_n}} \frac{\|p_nq\|_E}{\|p_nq\|_\Gamma}\right)^{1/k_n} = 1.
\end{equation}
Notice that for $\theta>0$ definitions (\ref{eq:ExPoly}) and (\ref{eq:ExPoly1}) coincide. Let
\begin{equation}
\nu^*(p_n)=\frac{1}{k_n} \sum_{p_n(z)=0} \delta_z
\label{10.29.07.3}
\end{equation}
and
\begin{equation}
\alpha^*(p_n)=\nu^*(p_n)_{|\overline D}+\widehat{\nu^*(p_n)_{|E\setminus \partial E}}.
\label{10.29.07.4}
\end{equation}
We remark that $ |\nu(p_n)|\le 1$ and $|\alpha^*(p_n)|\leq 1$. Let $\Nl(E) =\{ \nu:\nu \in \Lambda_1(E), \widehat \nu =\omega_{(E, \Gamma)}\},$ where  an $\omega_{(E,\Gamma)}$ is the Green equilibrium distribution on $E$ relative to $G$.

\begin{theorem}
\label{thm:2}
Let  $k_n\to\infty$ and $k_n=o(n)$ as $n\to\infty$. Then
\begin{equation}
\label{eq:L21}
\lim_{n\to\infty} \left(\frac{1}{k_n}\log\chi_n\right) = -\frac{1}{\cp(E,\Gamma)}.
\end{equation}
Moreover, if $\{p_n\}$ is an asymptotically extremal sequence  in the sense of (\ref{eq:ExPoly1}), then any weak-star limit point of $\{\nu^*(p_n)\}$ belongs to $\Nl(E)$ and 
\begin{equation}
\label{eq:L22}
\alpha^*(p_n) \cws \omega_{(E,\Gamma)} \quad \mbox{as} \quad n \to \infty.
\end{equation}
\end{theorem}

The last theorem provides the asymptotic behavior of the Kolmogorov $k$-width, $k=k_n$, of $A_n^\infty$ in $C(E)$. To formulate this theorem we need to introduce more notation. Fix $\theta\in(0,1]$ and define
\[
G_\theta := \left\{z\in\C:~U_D^{\lambda_\theta}(z)-g(z,\infty)>m_\theta\right\}.
\]
For $\theta=0$ we simply set $G_0:=G$. Clearly, the maximum principle for harmonic functions implies that $G_\theta=G$ whenever $S_\theta=\Gamma$ and it follows from (\ref{eq:m12}) and (\ref{eq:m10}) that $G\subseteq G_\theta$ for all $\theta\in[0,1]$. Let $G\subseteq G^\prime \subseteq G_\theta$, $H^\infty(G^\prime)$ be the space of bounded analytic functions on $G^\prime$, and $A_n^\infty(G^\prime)$ stand for the restriction to $E$ of the unit ball of $H^\infty(G^\prime)\cap\pl_n$. The following theorem shows that the $n$-th root limit of $d_{k_n}(A_n^\infty(G^\prime);C(E))$, $k_n/n\to\theta$, is independent of $G^\prime$.

\begin{theorem}
\label{thm:3}
Let $\{k_n\}_{n\in\N}$ satisfy (\ref{eq:theta}) for some $\theta\in[0,1]$, $G^\prime$ be a simply connected domain such that $G\subseteq G^\prime \subseteq G_\theta$, and $A_n^\infty=A_n^\infty(G^\prime)$. Then 
\begin{equation}
\label{eq:L31}
\lim_{n\to\infty} \left(\frac{1}{n}\log d_{k_n}(A_n^\infty;C(E))\right) = m_\theta.
\end{equation}
In particular, when $\theta=0$ and $k_n\to\infty$ as $n\to\infty$, we have that
\begin{equation}
\label{eq:L32}
\lim_{n\to\infty} \left(\frac{1}{k_n}\log d_{k_n}(A_n^\infty;C(E))\right) = -\frac{1}{\cp(E,\Gamma)}.
\end{equation}
\end{theorem}

\section{An Extremal Problem of the Potential Theory}
\label{sec:ContExtrPr}

We now state  the main  theorem of this section.  Let $\sigma$ be a compactly supported positive Borel  measure. Define
\[
M(\sigma) := \min_\Gamma U^\sigma-\min_E U^\sigma.
\]

\begin{theorem}
\label{thm:aux}
For each $\theta\in[0,1]$ we have
\begin{equation}
\label{eq:ExPr3}
m_\theta = \inf_{\mu\in\Lambda_\theta(E)}\sup_{\lambda\in\Lambda_{1-\theta}(\Gamma)} M(\mu+\lambda) = \sup_{\lambda\in\Lambda_{1-\theta}(\Gamma)}\inf_{\mu\in\Lambda_\theta(E)} M(\mu+\lambda).
\end{equation}
Moreover, if $\mu^*$, $|\mu^*|\leq\theta$, and $\lambda^*$, $|\lambda^*|\leq1-\theta$, are  compactly supported positive Borel measures such that
\begin{equation}
\label{eq:ExPr3U}
m_\theta = M(\mu^*+\lambda^*) = \sup_{\lambda\in\Lambda_{1-\theta}(\Gamma)} M(\mu^*+\lambda) 
\end{equation}
then $\supp(\mu^*)\subseteq E$, $\widehat{\mu^*}=\mu_\theta$, and $\lambda^*=\lambda_\theta$ when $S_\theta$ does not separate the plane and $\supp(\lambda^*)\subset\C\setminus G$ , $\widetilde{\lambda^*}=\lambda_\theta-(1-\theta -|\lambda^*|) \omega_\Gamma$ , otherwise.
\end{theorem}

The proof of Theorem $\ref {thm:aux}$ is based on several auxiliary lemmas.

\begin{lemma}
\label{lem:aux2}
We have
\begin{equation}
\label{eq:Lam0}
\lambda_0 = \omega_\Gamma, \;\;\;  \widehat\lambda_{0} = \omega_E, \;\;\; \mbox{and} \;\;\; m_0=0.
\end{equation}
\end{lemma}
\begin{proof} Since $U^{\omega_\Gamma}(z)=-\log\cp(\Gamma)$ for  $z\in\overline G$ and  $U^{\omega_E}(z)=-\log\cp(E)$ for $z\in E$,  it holds that
\begin{eqnarray}
U^{\omega_E}(z) = U^{\omega_\Gamma}(z)+c, \;\;\; z\in E, \nonumber
\end{eqnarray}
where $c=-\log\cp(E) + \log\cp(\Gamma)$. Using now the fact
$\supp (\omega_E) =\partial E$ and  the uniqueness of the balayage (see, for example, \cite[Thm. II.4.4]{SaffTotik}),  we can immediately conclude that 
\begin{equation}
\label{eq:21}
\widehat\omega_\Gamma = \omega_E 
\end{equation}
and 
\begin{equation}
\label{eq:200}
\int g(t,\infty)d\omega_\Gamma(t)=-\log\cp(E) + \log\cp(\Gamma).
\end{equation}
From this, on account of the formula
\begin{equation}
\label{eq:305}
U^{\omega_E}(z)=-\log\cp(E)-g(z,\infty), \;\;\; z\in \C,
\end{equation}
 we obtain that for every $z\in\C$,
\begin{eqnarray}
U_D^{\omega_\Gamma}(z)-g(z,\infty) &=& U^{\omega_\Gamma-\omega_E}(z)+\int g(t,\infty)d\omega_\Gamma(t)-g(z,\infty) \nonumber \\
\label{eq:22}
{} &=& U^{\omega_\Gamma}(z)+\log\cp(\Gamma).
\end{eqnarray}
So, 
\[
U_D^{\omega_\Gamma}(z)-g(z,\infty) =0, \;\;\; z\in \Gamma .
\]
Therefore,  relations $\lambda_0=\omega_\Gamma$ and $m_0=0$  follow from the uniqueness of the measure $\lambda_\theta$ satisfying conditions (\ref{eq:m12}) and (\ref{eq:m10})
(see Remark (e) after Theorem $\ref {thm:aux1}$).
\end{proof}

\begin{lemma}
\label{lem:aux3}
For each $\theta\in[0,1]$ we have
\begin{equation}
\label{eq:EqOnE}
\mu_\theta = \widehat{\lambda_0-\lambda_\theta} = \omega_E-\widehat{\lambda_\theta},
\end{equation}

\begin{equation}
\label{eq:34}
\widehat m_\theta=-\log\cp(E)-\int g(t, \infty) d \lambda_\theta(t),
\end{equation}
and
\begin{equation}
\label{eq:300}
U_D^{\lambda_\theta}(z) -g(z,\infty) = U^{\lambda_\theta+\mu_\theta}(z) - \widehat m_\theta, \;\;\; z\in\C.
\end{equation}
Moreover,
\begin{equation}
\label{eq:301}
\supp(\mu_\theta) = \partial E, \;\;\; \theta\in(0,1].
\end{equation}
\end{lemma}
\begin{proof} It is easy to see that for  $\theta=0$ and $\theta=1$ (\ref{eq:EqOnE}), (\ref{eq:34}) and (\ref{eq:300}) are valid.
In the case  $\theta=0$, $\mu_0=0$ by definition, and  Lemma $\ref{lem:aux3}$ implies that 
$\lambda_0 = \omega_\Gamma$ and $\widehat\lambda_0 = \omega_E$. From this, on account of $(\ref{eq:200})$ and $(\ref{eq:22})$ we get $(\ref{eq:34})$ and  $(\ref{eq:300})$.
For $\theta=1$,  $\lambda_1=0$ by definition.  It follows from Theorem $\ref{thm:aux2}$ (see Remark (a) after Theorem $\ref{thm:aux2}$)  that  $\mu_1=\omega_E$ and $\widehat m_1=-\log\cp(E)$. 
We also note that since  $\supp(\omega_E) = \partial E$, $(\ref{eq:301})$ holds  for $\theta=1$.

Let us consider now the case $\theta\in(0,1)$. We start from the next observation.  As noted in the Remark (b) after Theorem \ref{thm:1}, it  follows from \cite[Thm. 2.4]{LevLub01}, that $\lambda_0-\lambda_\theta$ is a positive measure. 
It is easy to see that $|\lambda_0-\lambda_\theta|=\theta.$ Hence, $\widehat{\lambda_0-\lambda_\theta}$ is a positive measure, and $|\widehat{\lambda_0-\lambda_\theta}|=\theta.$
Moreover, it is a simple application of the second unicity theorem \cite[Thm. II.4.6]{SaffTotik} to see that
\[
\widehat{\lambda_0 - \lambda_\theta} =\widehat\lambda_0 - \widehat\lambda_\theta .
\]
So, $\widehat\lambda_0 - \widehat\lambda_\theta =\omega_E-\widehat\lambda_\theta$ is a positive measure and $|\omega_E -\widehat{\lambda_\theta}|=\theta.$
According to the  property (\ref{10.29.07.1})  of the  Green potential, 
\begin{equation}
\label{eq:303}
U_D^{\lambda_\theta}(z) -g(z,\infty) = U^{\lambda_\theta-\widehat\lambda_\theta}(z)+\int g(t,\infty)d\lambda_\theta(t) - g(z,\infty),
\;\;\; z\in\C,
\end{equation}
and, by $(\ref{eq:305})$,
\begin{equation}
\label{eq:304}
U_D^{\lambda_\theta}(z) -g(z,\infty)=	U^{\lambda_\theta+\omega_E-\widehat\lambda_\theta}(z)+\log\cp(E)+\int g(t,\infty)d\lambda_\theta(t), \;\;\; z\in\C.
\end{equation}
Since
$U_D^{\lambda_\theta}(z)-g(z,\infty) =0$ on  $E$, (\ref {eq:EqOnE}) and (\ref {eq:34}) follow from the uniqueness  of the measure $\mu_\theta$ satisfying conditions (\ref{eq:m22}) and (\ref{eq:m20}). So, we have ($\ref{eq:300}$). Using now the facts that $E$ is a regular compact set, $\mu_\theta$ is the balayage of ${\lambda_0 - \lambda_\theta}$ and  properties of the balayage (see, for example,  \cite{Landkof}) we can conclude that  $\supp(\mu_\theta) = \partial E$.
\end{proof}

We can consider equation ($\ref{eq:300}$) as the basic equation of this section, it allows us to connect the Green potential $U_D^{\lambda_\theta}(z)-g(z,\infty)$ and the logarithmic potentials $U^{\lambda_\theta+\mu_\theta}(z)-	\widehat m_\theta$  of the extremal problems from Section \ref{sec:pot}. Since $U_D^{\lambda_\theta}(z)-g(z,\infty)=0$ on $E$ and 
$\min_{\Gamma}(U_D^{\lambda_\theta}(z)-g(z,\infty))= m_\theta$ on $\Gamma$,  we get immediately from ($\ref{eq:300}$) the  equality:
\begin{equation}
\label{eq:306}
m_\theta = M(\lambda_\theta+\mu_\theta).
\end{equation}

The function $U^{\lambda_\theta+\mu_\theta}$ satisfies the  following  property. The logarithmic potential $U^{\lambda_\theta+\mu_\theta}$   of a probability measure $\lambda_\theta+\mu_\theta$  is equal to constants on supports of $\mu_\theta$ and $\lambda_\theta$:
\begin{equation}
\label{10.25.07.1}
U^{\lambda_\theta+\mu_\theta} =	\widehat m_\theta \quad \mbox{on} \quad E \quad \mbox{and} \quad U^{\lambda_\theta+\mu_\theta} = \min_\Gamma U^{\lambda_\theta+\mu_\theta} = m_\theta+\widehat m_\theta \quad \mbox{on} \quad  S_{\theta}\subseteq \Gamma.
\end{equation}

\begin{lemma}
\label{lem:aux1}
For each  $ \theta \in [0, 1]$, we have
\begin{equation}
\label{eq:min2}
m_\theta = \inf_\mu M(\mu+\lambda_\theta), 
\end{equation}
where infimum is taken over all compactly supported positive Borel measures with $|\mu|\leq\theta$. Further, the equality in (\ref{eq:min2}), for $ \theta \in (0, 1]$,  is possible if and only if $\supp(\mu)\subseteq E$ and $\widehat\mu=\mu_\theta$.
\end{lemma}
\begin{proof} Let  $\theta=0$. In this case  $\mu_0=0$ by definition and by Lemma \ref{lem:aux2} $m_{0}=0$, $\lambda_0=\omega_\Gamma$. Since $U^{\omega_\Gamma}(z)=-\log\cp(\Gamma)$ for  $z\in\overline G$,  $ M(\lambda_0)=M(\omega_{\Gamma})=0$. This yields the equality $m_0=M(\lambda_0)$. 

Let  $\theta\in(0,1]$. Consider a logarithmic potential $U^{\mu-\mu_\theta}$. This function is superharmonic and bounded from below  in $D=\overline\C\setminus E$. Then by the generalized minimum principle for superharmonic functions \cite[Thm. I.2.4]{SaffTotik}, 
\begin{equation}
\label{eq:307}
\min_{E} U^{\mu-\mu_\theta} \le U^{\mu-\mu_\theta}(z), \;\;\;\;  z\in D.
\end{equation} 
In particular,
\begin{equation}
\min_E U^{\mu-\mu_\theta} \le \min_\Gamma U^{\mu-\mu_\theta}.
\label{10.25.07.2}
\end{equation}
Moreover, there are  strict inequalities in (\ref{eq:307}) for  $z \in D$ and in (\ref{10.25.07.2}), unless 
\[
\supp(\mu)\subseteq E ,\;\;\;  U^{\mu-\mu_\theta}(z)=0, \;\;\;  z\in D.
\]
That is if and only if $\widehat\mu=\mu_\theta$ by Carleson's unicity theorem (see \cite[Thm. II. 4.13]{SaffTotik}).

With the help of the equality 
\[
U^{\mu+\lambda_\theta}(z) = U^{\mu-\mu_\theta}(z) + U^{\lambda_\theta+\mu_\theta}(z), \quad z \in \C,
\]
and  (\ref{10.25.07.1}), we can write
\[
\min_E U^{\mu+\lambda_\theta} = \min_E U^{\mu-\mu_\theta} + \widehat m_\theta
\]
and
\[
\min_\Gamma U^{\mu+\lambda_\theta} \geq \min_\Gamma U^{\mu-\mu_\theta} + m_\theta +\widehat m_\theta.
\]
Therefore, by (\ref{10.25.07.2}),
\[
M(\mu+\lambda_\theta)=\min_\Gamma U^{\mu+\lambda_\theta}-\min_E U^{\mu+\lambda_\theta} \ge m_\theta,
\]
and the equality in (\ref{eq:min2}) is possible if and only if $\supp(\mu) \subseteq E$ and $\widehat \mu =\mu_\theta$.
\end{proof}

\begin{lemma}
\label{lem:aux4}
For each $\theta\in[0,1]$ we have
\begin{equation}
\label{eq:min3}
m_\theta = \sup_\lambda M(\mu_\theta+\lambda), 
\end{equation}
where supremum is taken over all compactly supported positive Borel measures with $ |\lambda|\leq1-\theta$. Further, the equality in (\ref{eq:min3}), for $\theta \in [0, 1)$ is possible if and only if $\lambda=\lambda_\theta$ when $S_\theta$ does not separate the plane and $\supp(\lambda)\subset\C\setminus G$, $\widetilde\lambda=\lambda_\theta-(1-\theta-|\lambda|) \omega_\Gamma$,  otherwise.
\end{lemma}
\begin{proof} Let $\theta=1$. In this case $\lambda_1=0$ and $m_1= -\max_\Gamma g(z, \infty)$
by definition and $\mu_1=\omega_E.$ On the basis of (\ref{eq:305}) we can write
\[
M(\mu_1)=\min_\Gamma U^{\omega_E}-\min_E U^{\omega_E}=-\max_{z \in \Gamma} g(z, \infty)=m_1.
\]
Let us consider the case when $\theta \in [0, 1).$
Denote by $\lambda$ any compactly supported positive Borel measure with mass at most $1-\theta$. It is enough to show that
\[
\min_{S_\theta}U^{\mu_\theta+\lambda} - \min_E U^{\mu_\theta+\lambda} \leq m_\theta.
\]

Consider a logarithmic potential $U^{\lambda-\lambda_\theta}$. This is superharmonic and bounded below function in $\overline\C\setminus S_\theta$. Then by the generalized minimum principle for superharmonic functions,
\begin{equation}
\min_{S_\theta} U^{\lambda-\lambda_\theta} \le U^{\lambda-\lambda_\theta}(z), \quad z \in \overline\C\setminus S_\theta, 
\label{10.25.07.3}
\end{equation}
and
\begin{equation}
\label{eq:40}
\min_{S_\theta} U^{\lambda-\lambda_\theta} \le \min_E U^{\lambda-\lambda_\theta}.
\end{equation}
From this, using  (\ref{10.25.07.1}),  we get 
\[
\min_{S_\theta} U^{\mu_\theta+\lambda}-\min_E U^{\mu_\theta+\lambda} \leq  m_\theta
\]
and therefore (\ref{eq:min3}) holds. Observe also that the equality in (\ref{eq:min3}) is possible if and only if we have the equality in (\ref{eq:40}). That is if and only if
\begin{equation}
U^{\lambda_\theta}(z) = U^{\lambda}(z) + c^*, \;\;\; z\in \Omega,
\label{11.06.07.1}
\end{equation}
where $c^*$ is some constant, $\Omega=\C\setminus S_\theta$ if $S_\theta$ does not separate the plane and $\Omega=G$ otherwise. In the former case $c^*=0$ and we get $\lambda=\lambda_\theta$ by Carleson's unicity theorem. In the latter situation $\supp (\lambda) \subseteq \C \setminus G$. Using the continuity of potentials in fine topology (see \cite[Sec. I.5]{SaffTotik}) and regularity of $\Gamma$, we may continue equality in (\ref{11.06.07.1}) up to $\overline G$. Let $\widetilde \lambda$ be the balayage of $\lambda$ onto $\Gamma$ relative to $\C\setminus G$ (we balayage only the part of $\lambda$ which is supported outside of $\overline G$). Then
\[
U^{\widetilde \lambda}(z)=U^\lambda(z)+c, \;\;\; z\in\overline G,
\]
where $\displaystyle c=\int g_{\overline G} (t, \infty) d \lambda(t).$  Thus,
\begin{equation}
\label{eq:41}
U^{\widetilde\lambda}(z)=U^{\lambda_\theta}(z)-c^* +c, \;\;\; z\in \overline G.
\end{equation}
Using now the maximum principle of harmonic functions in the domain $\overline \C \setminus \overline G$, we get
\begin{equation}
U^{\widetilde\lambda}(z) = U^{\lambda_\theta}(z) + \left(1-\theta-|\widetilde\lambda|\right)g_{\overline G}(z, \infty)-c^*+c, \;\;\; z\in\overline\C\setminus\overline G
\label{11.06.07.2}
\end{equation}
(we applied the maximum principle of harmonic functions for the difference of the left and right hand sides). Taking now on an account (\ref{eq:41}), we obtain the equality (\ref{11.06.07.2})  for all $z \in \overline \C$. From this with help of the  formula $\displaystyle U^{\omega_\Gamma}(z)=-\log \cp(\Gamma) -g_{\overline G}(z, \infty),$ and the unicity theorem \cite[Thm. II 2.1]{SaffTotik}, we can conclude that $\widetilde \lambda +(1-\theta -|\lambda|) \omega_\Gamma=\lambda_\theta$,
 which finishes the proof of the lemma.
\end{proof}
\begin{proof}[Proof of Theorem \ref{thm:aux}] It is a straightforward application of   Lemmas \ref{lem:aux1} and \ref{lem:aux4} to obtain
\[
\inf_{\mu\in\Lambda_\theta(E)}\sup_{\lambda\in\Lambda_{1-\theta}(\Gamma)} M(\mu+\lambda) \geq \inf_{\mu\in\Lambda_\theta(E)} M(\mu+\lambda_\theta) = m_\theta
\]
and
\[
\inf_{\mu\in\Lambda_\theta(E)}\sup_{\lambda\in\Lambda_{1-\theta}(\Gamma)} M(\mu+\lambda) \leq \sup_{\lambda\in\Lambda_{1-\theta}(\Gamma)} M(\mu_\theta+\lambda) = m_\theta.
\]
This establishes the first equality in (\ref{eq:ExPr3}). Clearly, we have
\[
\inf_{\mu\in\Lambda_\theta(E)}\sup_{\lambda\in\Lambda_{1-\theta}(\Gamma)} M(\mu+\lambda) \geq \sup_{\lambda\in\Lambda_{1-\theta}(\Gamma)}\inf_{\mu\in\Lambda_\theta(E)} M(\mu+\lambda).
\]
On the other hand, it follows from Lemmas \ref{lem:aux1} and \ref{lem:aux4} that
\[
\sup_{\lambda\in\Lambda_{1-\theta}(\Gamma)} M(\mu_\theta+\lambda) = \inf_{\mu\in\Lambda_\theta(E)} M(\mu+\lambda_\theta).
\]
Therefore, 
\[
\inf_{\mu\in\Lambda_\theta(E)}\sup_{\lambda\in\Lambda_{1-\theta}(\Gamma)} M(\mu+\lambda) \leq \sup_{\lambda\in\Lambda_{1-\theta}(\Gamma)}\inf_{\mu\in\Lambda_\theta(E)} M(\mu+\lambda),
\]
which finishes the proof of (\ref{eq:ExPr3}). Let now $\mu^*$ and $\lambda^*$ be as in (\ref{eq:ExPr3U}). Then by Lemma \ref{lem:aux1}, we observe that
\[
m_\theta  \leq M(\mu^*+\lambda_\theta) \leq \sup_\lambda M(\mu^*+\lambda)=m_\theta.
\]
Thus,
\[
m_\theta = M(\mu^*+\lambda_\theta)
\]
and $\supp (\mu^*)\subseteq E$ and $\widehat{\mu^*}=\mu_\theta$ again by Lemma \ref{lem:aux1}.  Furthermore, in this case
\[
m_\theta = M(\mu_\theta+\lambda^*)
\]
and, by Lemma \ref{lem:aux4}, $\lambda^*=\lambda_\theta$ when  $S_\theta$ does not separate the plane and  $\supp(\lambda^*)\subset \C\setminus G$, $\widetilde{\lambda^*}=\lambda_\theta-(1-\theta -|\lambda^*|) \omega_\Gamma$, otherwise.
\end{proof}

\section{Proofs of the Theorem \ref{thm:1}}
\label{sec:proofs1}

Before we present the proof of Theorem \ref{thm:1}, we introduce the analogue of the Tsuji points (\cite{Menk83}, \cite{FishS99}) that corresponds to the weighted  Green energy problem (\ref{eq:m11}). Set
\[
\delta_m^G := \max_{z_1,\ldots,z_m\in\Gamma}\left(\prod_{1\leq i<j\leq m}\exp\left\{-g(z_i,z_j)+\frac{g(z_i,\infty)}{1-\theta}+\frac{g(z_j,\infty)}{1-\theta}\right\}\right)^{2/m(m-1)}.
\]
Then
\[
\delta_m^G\geq\delta_{m+1}^G, \;\;\; m\in\N, \;\;\; \mbox{and} \;\;\; \lim_{m\to\infty}\log\delta_m^G = -J(\lambda_\theta)/(1-\theta)^2.
\]
Moreover, if $\{\zeta_1,\ldots,\zeta_m\}$ is any extremal set for $\delta_m^G$, then
\begin{equation}
\label{eq:m14}
\lambda_{m,\theta} \cws \lambda_\theta  \quad \mbox{as} \quad m \to \infty, \quad  \lambda_{m,\theta}:=\frac{1-\theta}{m}\sum_{j=1}^m\delta_{\zeta_j}.
\end{equation}
Here and in what follows we keep to the notation
\[
\lambda_n =\frac{1}{n} \sum_{j=1}^{n-k_n} \delta_{\xi_j, n-k_n}, 
\]
where $\{\xi_{1, n-k_n}, \dots, \xi_{n-n_k, n-k_n}\}$ is an extremal set for $\delta^G_{n-k_n}.$
We remark that
\begin{equation}
\lambda_n \cws \lambda_\theta \quad \mbox{as} \quad n \to \infty.
\label{10.30.07.1}
\end{equation}
The proof of these facts needs only minor modifications comparing to the case of the logarithmic kernel \cite[Thm. III.1.1-3]{SaffTotik}.

We also need a discretization of $\mu_\theta$. So, we introduce the Leja points (see \cite[Sec. III.1]{SaffTotik}) that correspond to the weighted minimal energy problem (\ref{eq:m21}). Set
\[
\delta_m := \max_{z_1,\ldots,z_m\in E}\left(\prod_{1\leq i<j\leq m}|z_i-z_j|\exp\left\{-\frac1\theta\left(U^{\lambda_\theta}(z_i)+U^{\lambda_\theta}(z_j)\right)\right\}\right)^{2/m(m-1)}.
\]
Then
\[
\delta_m\geq\delta_{m+1}, \;\;\; m\in\N, \;\;\; \mbox{and} \;\;\; \lim_{m\to\infty}\log\delta_m = -I(\mu_\theta)/\theta^2.
\]
Moreover, if $\{z_1, \ldots,z_m\}$ is any extremal set for $\delta_m$, then
\begin{equation}
\label{eq:m24}
\mu_{m,\theta} \cws \mu_\theta,  \;\;\; \mu_{m,\theta}:=\frac{\theta}{m}\sum_{j=1}^{m}\delta_{z_j}.
\end{equation}
Here and in what follows we keep to the notation
\[
\mu_n =\frac{1}{n} \sum_{j=1}^{k_n} \delta_{z_{j, k_n}}, 
\]
where $\{z_{1, k_n}, \dots, z_{k_n, k_n}\}$ is an external set for $\delta_{k_n}$. We have
\begin{equation}
\mu_n \cws \mu_\theta \quad \mbox{as} \quad n \to \infty.
\label{10.30.07.2}
\end{equation}

It is easy to see that for any compact set $K$, $p\in\pl_{k_n}$, $p\not \equiv 0$, and $q\in\pl_{n-k_n}$, $q\not \equiv 0$,  we have
\[
\|pq\|_K^{1/n}=\gamma^{1/n}\exp\left\{-\min_KU^{\nu(p)+\nu(q)}\right\},
\]
where $\gamma$ is the leading coefficient of $pq$,
and 
\[
\nu(p):=\frac1n\sum_{p(z)=0}\delta_z \;\;\; \mbox{and} \;\;\; \nu(q):=\frac1n\sum_{q(z)=0}\delta_z
\]
(the sums are taken counting multiplicities of zeros of $p$ and $q$). Therefore, we get
\[
\frac1n\log\left(\frac{\|pq\|_E}{\|pq\|_\Gamma}\right) =  \min_\Gamma U^{\nu(p)+\nu(q)}-\min_E U^{\nu(p)+\nu(q)} = M(\nu(p)+\nu(q)).
\]
\begin{proof}[Proof of Theorem \ref{thm:1}] Let 
\[
\sup_{q \in \pl_{n-k_n}} M(\mu_n+\nu(q))=M(\mu_n+\nu(Q_n))
\]
for some polynomial $Q_n \in \pl_{n-k_n}, Q_n \not \equiv 0.$ Denote by
\[
\sigma_n=\nu(Q_n)_{|\overline G}+\widetilde{\nu(Q_n)_{|\C\setminus \overline G}}.
\]
By properties of the balayage, $\supp (\sigma_n) \subseteq \overline G, \ |\sigma_n|=|\nu(Q_n)|, $ and
\[
M(\mu_n+\nu(Q_n))=M(\mu_n+\sigma_n).
\]
We now choose  a convergent subsequence such that  
\begin{equation}
\sigma_n \cws \sigma, \quad n \in \Lambda \subset \N,
\label{10.21.07.1}
\end{equation}
and
\[
\limsup_{n \to \infty} M(\mu_n+\sigma_n)=\lim_{n \to \infty,  n \in \Lambda} M(\mu_n+\sigma_n).
\]
We remark that $\supp(\sigma) \subseteq \overline G$ and $|\sigma| \le 1-\theta$.
Since $E$ and $\Gamma$ are regular sets, conditions  (\ref{10.30.07.2}) and (\ref{10.21.07.1})  imply (cf. \cite{GR1}) that
\[
\min_\Gamma  V^{\mu_n+\sigma_n} \to \min_\Gamma V^{\mu_\theta+\sigma} \quad  \mbox {as} \quad n \to \infty,\quad n \in \Lambda,
\]
\[
\min_E  V^{\mu_n+\sigma_n} \to \min_E V^{\mu_\theta+\sigma} \quad  \mbox {as} \quad n \to \infty, \quad n \in \Lambda,
\]
and then
\[
\lim_{n \to \infty, n \in \Lambda} M(\mu_n+\sigma_n) =M(\mu_\theta+\sigma), 
\]
where, by Lemma \ref{lem:aux4}, $M(\mu_\theta+\sigma) \le m_\theta$. Therefore,
\begin{eqnarray}
\limsup_{n\to\infty} \left(\frac{1}{n} \log \chi_n\right) &\leq& \limsup_{n\to\infty} \sup_{q\in\pl_{n-k_n}} M(\mu_n+\nu(q)) =\limsup_{n \to \infty} M(\mu_n+\nu(Q_n))\nonumber\\
&=& \lim_{n\to\infty, n \in \Lambda}  M(\mu_\theta+\sigma_n)  =
 M(\mu_0+\sigma) \le  m_\theta.
\label{eq:t11}
\end{eqnarray}
For any polynomial $p \in \pl_{k_n}, p \not \equiv 0, $ consider the following function
\[
u(z)=U^{\nu(p)+\lambda_n}(z) -\min_E U^{\nu(p)+\lambda_n} -U_D^{\lambda_n}(z) +g(z, \infty), \quad z \in \C.
\]
This function is superharmonic in $D$. Using the generalized minimum principle for superharmonic functions, we obtain that $u(z) \ge 0$, $z \in D$. In particular,
\begin{equation}
\min_\Gamma U^{\nu(p) +\lambda_n} -\min_E U^{\nu(p)+\lambda_n} \ge \min_\Gamma \left(U^{\lambda_n}_D(z) -g(z, \infty)\right).
\label{11.01.07.1}
\end{equation}
Since (\ref{11.01.07.1}) is valid for any $p \in \pl_{k_n}, p \not \equiv 0, $ we get
\begin{equation}
\frac{1}{n} \log \chi_n \ge \min_{\Gamma} \left(U^{\lambda_n}_D(z) -g(z, \infty)\right).
\label{11.01.07.2}
\end{equation}
Further, in view of the properties of weakly convergent sequences,
\[
\min_\Gamma \left(U^{\lambda_n}_D(z)-g(z, \infty)\right) \to \min_\Gamma \left(U^{\lambda_\theta}_D(z) -g(z, \infty)\right)=m_\theta \quad \mbox{as} \quad n \to \infty.
\]
Then, by the relation (\ref{11.01.07.2}), we get
\begin{equation}
\liminf_{n \to \infty} \left(\frac{1}{n} \log \chi_n\right) \ge m_\theta.
\label{eq:t12}
\end{equation}
So, (\ref{eq:L11}) follows from (\ref{eq:t11}) and  (\ref{eq:t12}).
 
Fix a positive $R$ such that $\overline G\subset U$, where $U=\{z:|z| <R\}$. Let $L=\{z:|z|=R\}$.

Let now $\{p_n,q_n\}_{n\in\N}$, $p_n\in\pl_{k_n}$ and $q_n\in\pl_{n-k_n}$, be a sequence of asymptotically extremal pairs of polynomials. First, we show (\ref{eq:L12}). Let $\{\nu(p_n)\}$, $n \in \Lambda_0 \subset \N$, be a convergent subsequence. Let 
\begin{equation}
\sigma_n =\nu(p_n)_{|\overline U}+\tau_n, 
\label{11.03.07.4}
\end{equation}
where $\tau_n$ is the balayage of $\nu(p_n)_{|\C\setminus \overline U}$ on $L$. By  the properties  of balayage,  for any polynomial $q \in \pl_{n -k_n}$, $ q \not \equiv 0, $
\[
M(\nu(p_n)+\nu(q)) =  M(\sigma_n+\nu(q)).
\]
Hence,
\begin{equation}
\sup_{q \in \pl_{n-k_n}} M(\nu(p_n)+\nu(q)) \ge M(\sigma_n+\lambda_n).
\label{11.01.07.3}
\end{equation}
We choose  a convergent subsequence
\begin{equation}
\sigma_n \stackrel{*}{\to} \nu, \;\; n \in \Lambda \subset \Lambda_0\subset \N, 
\label{11.03.07.5}
\end{equation}
where  $|\nu| \le \theta, \supp (\nu) \subseteq \overline U.$ On the basis of the  fact that $E$ and $\Gamma $ are regular sets, we get
\begin{equation}
M(\sigma_n+\lambda_n) \to M(\nu+\lambda_\theta) \quad \mbox{as} \quad n \to \infty, \quad n \in \Lambda.
\label{11.01.07.4}
\end{equation}
By (\ref{eq:ExPoly}), 
\[
\limsup_{n \to \infty} \sup_{q \in \pl_{n-k_n}} M((p_n)+\nu(q))=m_\theta.
\]
From this, on an account of (\ref{11.01.07.3}) and (\ref{11.01.07.4}), we obtain that $M(\nu+\lambda_\theta)\le m_\theta$. Applying  Lemma \ref{lem:aux1}, we can write $M(\nu+\lambda_\theta)=m_\theta$, $\supp (\nu) \subseteq E$, and $\widehat \nu=\mu_\theta$.

Since $\supp (\nu) \subseteq E$, we obtain from (\ref{11.03.07.4}) and (\ref{11.03.07.5}) that
\[
\tau_n \cws 0 \quad \mbox{as} \quad n \to \infty, \quad n \in \Lambda,
\]
\[
\nu(p_n) \cws \nu \quad \mbox{as} \quad n \to \infty, \quad n \in \Lambda,
\]
and, then,
\[
\nu(p_n) \cws \nu \quad \mbox{as} \quad n \to \infty, \quad n \in \Lambda_0.
\]
From this, by properties of the balayage,
\[
\alpha(p_n) \cws \mu_\theta \quad \mbox{as} \quad n \to \infty,\quad n \in \Lambda_0,
\]
and, then, 
\[
\alpha(p_n) \cws \mu_\theta \quad \mbox{as} \quad n \to \infty.
\]
The relation (\ref{eq:L12}) thereby is obtained.

It only remains to prove (\ref{eq:L13}). By properties of the balayage, $M(\nu(p_n)+\nu(q_n))=M(\alpha(p_n)+\nu(q_n)).$ Define 
\begin{equation}
\nu_n=\nu(q_n)_{|\overline U}+\eta_n,
\label{11.03.07.6}
\end{equation}
where $\eta_n$ is the balayage of $\displaystyle \nu(q_n)_{|\C\setminus \overline U}$ on $L$, when $S_\theta$ does not separate the plane and $\nu_n=\beta(q_n)$  otherwise. Then  $M(\alpha(p_n)+\nu(p_n))=M(\alpha(p_n)+\nu_n)$ and  (\ref{eq:ExPoly}) yields that
\[
\lim_{n\to\infty} M(\alpha(p_n)+\nu_n) = m_\theta.
\]
As above, taking a convergent subsequence, $\nu_n \stackrel{*}{\to} \nu$, $n \in \Lambda \subset \N$, we get
\[
\lim_{n \to \infty, n \in \Lambda} M(\alpha(p_n)+\nu_n)=M(\mu_\theta+\nu)
\]
and $M(\mu_\theta+\nu)=m_\theta.$

Let us consider now the case when $S_\theta$ does not separate the plane. Since, by Lemma \ref{lem:aux4}, $ \nu =\lambda_\theta$, $\supp (\lambda_\theta) \subseteq \Gamma$, we obtain that
\[
\eta_n \cws 0 \quad \mbox{as} \quad n \to \infty, \quad n \in \Lambda, 
\]
\[
\nu(q_n) \cws \nu=\lambda_\theta \quad \mbox{as} \quad n \to \infty, \quad n \in \Lambda,
\]
and then
\[
\nu(q_n) \cws \lambda_\theta \quad \mbox{as} \quad n \to \infty.
\]

In the case when $S_\theta$ does separate the plane, we have by Lemma \ref{lem:aux4} that $\supp (\nu) \subset \C\setminus G$, $|\nu|=1$, and  $\widetilde \nu =\lambda_\theta$. Moreover, by definition  of $\beta (q_n)$ (see (\ref{11.03.07.2}))  we can conclude that $\supp (\nu) \subseteq \Gamma$. From this  and the fact that $\widetilde \nu =\lambda_\theta$ we  obtain that $\nu=\lambda_\theta$ and then
\[
\beta(q_n) \stackrel{*}{\to}  \lambda_\theta \quad \mbox{as} \quad n \to \infty. \qedhere
\]
\end{proof}

\section{Some properties of $\lambda_\theta$ and $\mu_\theta$}
\label{sec:lm}

In the next two theorems we describe some properties of the extremal  measures $\lambda_\theta$ and $\mu_\theta$, their supports, and the constants $m_\theta$ and $\widehat m_\theta$. It will be convenient for us to use the notation
\[
\omega_{(K,\partial E)} \;\;\; \mbox{and} \;\;\; \cp(K,\partial E)
\]
for the Green equilibrium distribution and the condenser capacity of a compact set $K\subset D$, respectively (cf. \cite[Ch. II and VII]{SaffTotik}).

\begin{theorem}
\label{thm:4}
\begin{itemize}
\item[(a)]  The family $\{S_\theta\}, \ \theta\in [0, 1],$ is a decreasing family of sets such that
\[
S_\theta = \overline{\bigcup_{\theta<\tau<1}S_\tau} \subseteq \bigcap_{0\leq\tau<\theta}S_\tau =  \left\{z\in\Gamma:~U_D^{\lambda_\theta}(z)-g(z,\infty)=m_\theta\right\}
\]
and
\[
S_1 := \bigcap_{0\leq\tau<1}S_\tau = \left\{z\in\Gamma:~g(z,\infty)=-m_1\right\};
\]
\item [(b)] the family $\{\lambda_\theta\}, \ \theta\in [0, 1), $ is decreasing and continuous in the weak$^*$ sense. Moreover,
\[
\lambda_\theta = \int_\theta^1\omega_{(S_\tau,\partial E)}d\tau;
\]
\item [(c)] $m_\theta$ is a continuous and strictly decreasing function of $\theta$ on $[0,1]$. Furthermore,
\[
m_\theta = m_1+\int_\theta^1\frac{d\tau}{\cp(S_\tau,\partial E)}
\]
and
\[
\frac{m_\theta}{\theta}\to - \frac{1}{\cp(\Gamma,\partial E)} \;\;\; \mbox{as} \;\;\; \theta\to0;
\]
\item[(d)] \quad  $\displaystyle \frac{\lambda_0-\lambda_\theta}{\theta} \cws \omega_{(\Gamma, \partial E)}$ as  $\theta \to 0$;
\item[(e)] if $S_1$ has positive capacity then $\displaystyle \frac{\lambda_\theta}{1-\theta}\stackrel{*}{\to} \omega_{(S_1,\partial E)}$ as $\theta\to1$.
\end{itemize}
\end{theorem}
\begin{proof} Statements (a), (b), and the first part of (c) follow from \cite[Thm. 2.4]{LevLub01}. (We should remark that it is required in \cite[Thm. 2.4]{LevLub01} for any compact set $K\subset\Gamma$ to have connected complement. However, a direct examination of the proof shows that the theorem still holds when $E$ is contained in $G$ and the later is simply connected.) Further, by Lemma \ref{lem:aux2}, $m_0=0$ and $S_0=\Gamma$. This means that
\[
\frac{m_\theta}{\theta}=-\frac1\theta\int_0^\theta\frac{d\tau}{\cp(S_\tau,\partial E)}.
\]
Thus, the second part of (c) follows by the continuity of $\cp(S_\theta,\partial E)$ as function of $\theta$ at zero from the right \cite[Thm. 5.1.3]{Ransford}.

It has been proved in \cite[Thm. 2.4]{LevLub01} that
\begin{equation}
\frac{d\lambda_\theta}{d\theta}=-\omega_{(S_\theta,\partial E)}
\label{10.25.07.4}
\end{equation}
for any point of continuity of $\cp(S_\theta,\partial E)$ as a function of $\theta$. Then (d) follows from continuity of $\cp(S_\theta, \partial E)$ at $\theta=0$ and the fact $S_0=\Gamma.$ Now, assume that $S_1$ has positive logarithmic capacity and therefore well-defined Green equilibrium distribution $\omega_{(S_1,\partial E)}$. As above, we can use (\ref{10.25.07.4}). The continuity from the left of $\cp(S_\theta,\partial E)$ at one follows from \cite[Thm. 5.1.3]{Ransford} by the  definition of $S_1$.
\end{proof}

The following theorem describes  the connection between $\mu_\theta$ and $\lambda_\theta$, some properties of $\mu_\theta$ and $ \widehat{m_\theta}$.

\begin{theorem}
\label{thm:5}
\begin{itemize}
\item [(a)] The family $\{\mu_\theta\}, \ \theta\in(0,1], $ is increasing, continuous in the weak$^*$ sense, and such that
\[
\supp(\mu_\theta)=\partial E \;\;\; \mbox{and} \;\;\; \mu_\theta = \omega_E - \widehat\lambda_\theta.
\]
Moreover,
\[
\frac{\mu_\theta}{\theta}\stackrel{*}{\to}\omega_{(E,\Gamma)} \;\;\; \mbox{as} \;\;\; \theta\to0;
\]
\item [(b)] $\widehat m_\theta$ is a continuous and strictly increasing function of $\theta$ on $[0,1]$. Furthermore,
\[
\widehat m_\theta = -\log\cp(E)-\int g(t, \infty) d \lambda_\theta(t).
\]
\end{itemize}
\end{theorem}
\begin{proof} Part (b) follows from (\ref{eq:34}) and Theorem \ref{thm:4}(b). First part of (a) follows from  Lemma \ref{lem:aux3}, the formula $\mu_\theta=\widehat{\lambda_0-\lambda_\theta}$ (see (\ref{eq:EqOnE})), and Theorem \ref{thm:4}(b). The continuity of $\{\mu_\theta\}$  follows from  continuity of the family $\{\lambda_\theta\}$, formula $\mu_\theta =\widehat{\lambda_0-\lambda_\theta}$, and properties of the balayage (see, for example \cite{Landkof}).
We have
\[
\frac{\lambda_0-\lambda_\theta}{\theta} \stackrel{*}{\to} \omega_{(\Gamma,\partial E)} \;\;\; \mbox{as} \;\;\; \theta\to0.
\]
Thus, by properties of balayage, 
\[
\widehat{\frac{\lambda_0-\lambda_\theta}{\theta} } \stackrel{*}{\to} \widehat{\omega_{(\Gamma, \partial E)}} \quad \mbox{as} \quad \theta \to 0, 
\]
and then
\[
\mu_\theta/\theta \stackrel{*}{\to} \widehat{\omega_{(\Gamma,\partial E)}} \;\;\; \mbox{as} \;\;\; \theta\to0.
\]
It remains only to remark that
\[
\widehat{\omega_{(\Gamma,\partial E)}} = \omega_{(E,\Gamma)}. \qedhere
\]
\end{proof}

\section{Proof of Theorem \ref{thm:2}}
\label{sec:proofs2}

\begin{proof}[Proof of Theorem \ref{thm:2}] Since $k_n=o(n)$ as $ n \to \infty$, $[k_n/\theta]\leq n$ for any fixed $\theta\in(0,1)$ and  $n$ sufficiently large. Let $l_n=[k_n/\theta].$ Therefore,
\[
\inf_{p \in \pl_{k_n}} \sup_{q \in \pl_{l_n-k_n}} \frac{\|pq\|_E}{\|pq\|_\Gamma} \leq \chi_n \leq d_{k_n}(A^\infty;C(E))
\]
(compare (\ref{eq:kWidthHol}) with the definition (\ref{eq:ExPr1}) of $\chi_n$). Then by (\ref{eq:L11}) and (\ref{eq:kWidthLim}), we have
\begin{equation}
\label{eq:t20}
\frac{m_\theta}{\theta} \leq  \liminf_{n\to\infty}\left(\frac{1}{k_n}\log\chi_n\right)\le \limsup_{ n \to \infty} \left(\frac{1}{k_n} \log \chi_n\right) \leq-\frac{1}{\cp(E,\Gamma)}.
\end{equation}
Taking the limit $\theta \to 0$, we obtain (\ref{eq:L21}) from Theorem \ref{thm:4}(c) and the fact that $\cp(E,\Gamma)=\cp(\Gamma,\partial E)$.

Let now $p_n$ be asymptotically extremal polynomials in the sense of (\ref{eq:ExPoly1}). Fix an arbitrary $\theta\in(0,1)$. Let
\[
\sigma_n=\frac{1}{l_n} \sum_{j=1}^{l_n-k_n} \delta_{\xi_j, l_n-k_n}, 
\]
where $\{\xi_{1, l_n-k_n}, \dots, \xi_{l_n-k_n, l_n-k_n}\}$ is an extremal set for $\delta_{l_n-k_n}^G$. Observe that
\[
\sigma_n \cws \lambda_\theta \quad \mbox{as} \quad n \to \infty.
\]
Denote by $\displaystyle Q_n(z) =\prod_{j=1}^{l_n-k_n} (z-\xi_{j, l_n-k_n})$ the corresponding  polynomial degree $l_n-k_n$.  By (\ref{eq:ExPoly1}),
\begin{equation}
\limsup_{n \to \infty} \left( \frac{1}{\chi_n} \frac{||p_nQ_n||_E}{||p_nQ_n||_\Gamma}\right)^{1/k_n} \le 1.
\label{11.01.07.5}
\end{equation}
We have
\begin{equation}
\frac{1}{k_n} \log\frac{||p_nQ_n||_E}{||p_nQ_n||_\Gamma}=M\left(\nu^*(p_n)+\frac{l_n}{k_n}\sigma_n\right).
 \label{11.01.07.6}
\end{equation}
As in the proof of Theorem \ref{thm:1}, we fix $R>0$ such that $\overline G \subset U =\{z:|z|<R\}$. Let $L=\{z:|z|=R\}$. Let
\[
\nu_n=\nu^* (p_n)_{|\overline U}+\eta_n, 
\]
where $\eta_n$ is the balayage of $\nu^*(p_n)_{|\C\setminus \overline U}$ onto $L$. According to the properties of the balayage,
\[
M\left(\nu^*(p_n)+\frac{l_n}{k_n} \sigma_n\right)=M\left(\nu_n+\frac{l_n}{k_n} \sigma_n\right).
\]
By (\ref{eq:L21}) and (\ref{11.01.07.5}), 
\begin{equation}
\label{eq:t21}
\limsup_{n\to\infty} M\left(\nu_n+\frac{l_n}{k_n}\sigma_n\right) \leq -\frac{1}{\cp(E,\Gamma)}.
\end{equation}
We select a convergent subsequence $\nu_n \cws \nu$, $n \in \Lambda \subset \N$, $|\nu|\le 1$, such that 
\[
\limsup_{n \to \infty} M\left(\nu_n+\frac{l_n}{k_n}\sigma_n\right)=\lim_{n \to \infty, n \in \Lambda} M\left(\nu_n+\frac{l_n}{k_n}\sigma_n\right).
\]
Since $E$ and $\Gamma$ are regular sets,
\[
\lim_{n \to \infty, n \in \Lambda}  M\left(\nu_n+\frac{l_n}{k_n}\sigma_n\right)= M \left(\nu+\frac{\lambda_\theta}{\theta}\right).
\]
Therefore, by (\ref{eq:t21}), we get
\begin{equation}
\label{eq:t22}
 M\left(\nu+\frac1\theta\lambda_\theta\right) \leq -\frac{1}{\cp(E,\Gamma)}.
\end{equation}
Since $\lambda_0=\omega_\Gamma$ and $\displaystyle U^{\omega_\Gamma}(z)=-\log \cp (\Gamma)$ on $\overline G$, 
\[
M\left( \nu +\frac{1}{\theta}\lambda_\theta \right) =M\left(\nu -\frac{\lambda_0-\lambda_\theta}{\theta}\right)
\]
and
\[
M\left(\nu -\frac{\lambda_0-\lambda_\theta}{\theta}\right)\le -\frac{1}{\cp(E, \Gamma)}.
\]
According to Theorem \ref{thm:4}(d),
\[
\frac{\lambda_0-\lambda_\theta}{\theta} \stackrel{*}{\to} \omega_{(\Gamma, \partial E)} \quad \mbox{as} \quad \theta \to 0. 
\]
Taking now the limit as $ \theta \to 0$, we get
\begin{equation}
M(\nu-\omega_{(\Gamma, \partial E)}) \le -\frac{1}{\cp(E, \Gamma)}.
\label{11.07.07.1}
\end{equation}
Since
\[
U^{\omega_{(E, \Gamma)}-\omega_{(\Gamma, \partial E)}}=\frac{1}{\cp(E,\Gamma)} \quad \mbox{on} \quad E
\]
and
\[
U^{\omega_{(E, \Gamma)}-\omega_{(\Gamma, \partial E)}}=0 \quad \mbox{on} \quad \Gamma,
\]
we obtain that
\[
M(\nu-\omega_{(\Gamma, \partial E)})=M(\nu-\omega_{(E, \Gamma)})-\frac{1}{\cp(E, \Gamma)}.
\]
Thus, we derive from this and (\ref{11.07.07.1}) that
\[
M(\nu-\omega_{(E, \Gamma)}) \le 0
\]
and
\begin{equation}
\label{eq:t24}
\min_\Gamma U^{\nu-\omega_{(E,\Gamma)}}\le \min_E U^{\nu-\omega_{(E,\Gamma)}}.
\end{equation}
Applying now the generalized minimum principle for superharmonic functions,  we can conclude that 
$\min_\Gamma U^{\nu-\omega_{(E,\Gamma)}}=\min_E U^{\nu-\omega_{(E,\Gamma)}},$
$U^{\nu-\omega_{(E, \Gamma)}}=0$ in $D$ and $\supp (\nu) \subseteq E.$ By Carleson's unicity theorem, we obtain from this that $\widehat \nu=\omega_{(E, \Gamma)}$. Since $\supp (\nu) \subseteq E$, we get
\[
\eta_n \cws 0 \quad \mbox{as} \quad n \to \infty, \quad n \in \Lambda,
\]
and
\[
\nu^*(p_n) \cws \nu \quad \mbox{as} \quad n \to \infty, \quad n \in \Lambda.
\]
From this, by the properties of the balayage, we can write
\[
\alpha^*(p_n) \cws \omega_{(E, \Gamma)} \quad \mbox{as} \quad n \to \infty, \quad n \in \Lambda,
\]
and then
\[
\alpha^*(p_n) \cws \omega_{(E, \Gamma)} \quad \mbox{as} \quad n \to \infty. \qedhere
\]
\end{proof}

\section{Proof of Theorem \ref{thm:3}}
\label{sec:proofs3}

\begin{proof}[Proof of Theorem \ref{thm:3}] We start by showing the lower bounds in (\ref{eq:L31}) and (\ref{eq:L32}). Since 
\[
A_n^\infty(G_\theta) \subseteq A_n^\infty(G^\prime), \quad G^\prime \subseteq G_\theta,
\]
we may take $G^\prime=G_\theta$.

Let $\{\lambda_n\}$ be a sequence of measures defined as in (\ref{10.30.07.1}). For each $\theta\in[0,1)$ we take $\{q_n\}$ to be the sequence of monic polynomials such that $\nu(q_n)=\lambda_n$, where $\nu(h)$ be the counting measure of the zeros of a polynomial $h$ normalized by $1/n$. For $\theta=1$ we take $q_n\equiv1$. Then $\{pq_n:~p\in\pl_{k_n}\}$ is a linear space of continuous functions on $E$ of dimesion $k_n+1$. Hence, it follows from \cite[pg. 137]{Lorentz} that for any linear space of continuous functions on $E$ of dimension $k_n$, say $X_{k_n}$, there exists a polynomials $p_{X_{k_n}}$ such that
\[
\inf_{g \in X_{k_n}}\|p_{X_{k_n}}q_n-g\|_E \geq \|p_{X_{k_n}}q_n\|_E.
\]
In particular, it means that
\begin{equation}
\label{eq:t39}
d_{k_n}(A_n^\infty;C(E)) \geq \inf_{X_{k_n}} \frac{\|p_{X_{k_n}}q_n\|_E}{\|p_{X_{k_n}}q_n\|_{\Gamma_\theta}} \geq \inf_{p\in\pl_{k_n}}\frac{\|pq_n\|_E}{\|pq_n\|_{\Gamma_\theta}}.
\end{equation}
When $\theta=1$, we get from the Bernstein-Walsh inequality and (\ref{eq:t39}) that
\begin{equation}
\label{eq:t36}
d_{k_n}(A_n^\infty;C(E)) \geq \exp\{k_nm_1\}.
\end{equation}
For $\theta\in[0,1)$, the lower estimate in (\ref{eq:t39}) yields
\begin{eqnarray}
d_{k_n}(A_n^\infty;C(E)) & \geq & \inf_{p\in\pl_{k_n}} \exp\left\{n\left(\min_{\Gamma_\theta}U^{\nu(p)+\lambda_n}-\min_EU^{\nu(p)+\lambda_n}\right)\right\} \nonumber \\ \nonumber \\
\label{eq:t37}
{} &\geq& \exp\left\{n\min_{\Gamma_\theta} \left(U_D^{\lambda_n}(z)-g(z,\infty)\right)\right\},
\end{eqnarray}
where $\Gamma_\theta := \partial G_\theta$ and we used (\ref{11.01.07.1}) with $\Gamma_\theta$ instead of $\Gamma$. As before, by the properties of weakly convergent sequences, it holds that
\begin{equation}
\label{eq:t312}
\min_{\Gamma_\theta} \left(U^{\lambda_n}_D(z)-g(z, \infty)\right) \to \min_{\Gamma_\theta} \left(U^{\lambda_\theta}_D(z) -g(z, \infty)\right)=m_\theta \quad \mbox{as} \quad n \to \infty.
\end{equation}
Thus, we get from (\ref{eq:t36}) and (\ref{eq:t37}) with (\ref{eq:t312}) that
\begin{equation}
\label{eq:t38}
\liminf_{n\to\infty}\left(\frac1n\log d_{k_n}(A_n^\infty;C(E))\right) \geq m_\theta.
\end{equation}
When $\theta=0$, we have that $G_0=G$. Further, we get exactly as in the first inequality in (\ref{eq:t20}) that
\[
\liminf_{n\to\infty}\left(\frac{1}{k_n}\log d_{k_n}(A_n^\infty;C(E))\right) \geq \liminf_{n\to\infty}\left(\frac{l_n}{k_n}\frac{1}{l_n}\log d_{k_n}(A_{l_n}^\infty(G_\tau);C(E))\right) \geq \frac{m_\tau}{\tau}
\]
for any $\tau\in(0,1]$, where $l_n:=[k_n/\tau]$ and we used (\ref{eq:t38}) and the fact that $G\subseteq G_\tau$. Therefore,
\begin{equation}
\label{eq:t310}
\liminf_{n\to\infty}\left(\frac{1}{k_n}\log d_{k_n}(A_n^\infty;C(E))\right) \geq -\frac{1}{\cp(E,\Gamma)}
\end{equation}
by Theorem \ref{thm:4}(c).

Now we shall show the upper bounds in (\ref{eq:L31}) and (\ref{eq:L32}). Observe that
\[
d_{k_n}(A_n^\infty;C(E)) \leq d_{k_n}(A^\infty;C(E)).
\]
Thus, (\ref{eq:L32}) follows from (\ref{eq:t310}) and (\ref{eq:kWidthLim}). Since for $\theta=0$ limit (\ref{eq:L31}) follows from (\ref{eq:L32}), we may assume that $\theta\in(0,1]$. Moreover, since 
\[
A_n^\infty(G^\prime) \subseteq A_n^\infty(G), \quad G \subseteq G^\prime,
\]
we may take $G^\prime=G$. To proceed with the upper bound we need to construct a special sequence of domains. Fix $\theta\in(0,1]$ and define
\[
\Omega_{\theta,\delta} := \left\{z\in\C:~U_D^{\lambda_\theta}(z)-g(z,\infty) < m_\theta+\delta\right\}, \quad \delta\in(0,-m_\theta).
\]
Each such domain $\Omega_{\theta,\delta}$ is unbounded and contains $S_\theta=\supp(\lambda_\theta)$ by (\ref{eq:m12}). Also denote
\[
G^\delta := \left\{z\in\C:~g_{\overline G}(z,\infty)\leq\delta\right\},
\]
where $g_{\overline G}(\cdot,\infty)$ is the Green function with pole at infinity for $\overline\C\setminus G$. Now, for each fixed $\delta\in(0,-m_\theta)$ take $U_\delta$ to be a connected domain (possibly unbounded) with regular boundary and such that
\[
S_\theta \subset \overline\C\setminus\overline U_\delta, \quad E\subset U_\delta, \quad \mbox{and} \quad L_\delta:=\partial U_\delta \subset \Omega_{\theta,\delta}\cap G^\delta.
\]
Then the harmonic measure (cf. \cite[Sec. 4.3]{Ransford}) for $U_\delta$, say $\omega_\delta(\cdot,\cdot)$, exists,
\begin{equation}
\label{eq:t34}
\|h\|_{L_\delta} \leq \|h\|_{\partial G^\delta} \leq \|h\|_\Gamma\exp\{n\delta\}, \quad h\in\pl_n,
\end{equation}
by the Bernstein-Walsh inequality, and
\begin{equation}
\label{eq:t35}
\max_{L_\delta}\left(U_D^{\lambda_\theta}(z)-g(z,\infty)\right) \leq m_\theta+\delta.
\end{equation}
Therefore, if $S_\theta=\Gamma$, then $U_\delta$ is an open subset of $G$ that contains $E$ and whose boundary  is regular and close enough to $\Gamma$ so (\ref{eq:t35}) holds. If $S_\theta$ is a proper subset of $\Gamma$, the $U_\delta$ is an unbounded open set that contains $E$, whose boundary is regular, encompasses $S_\theta$, and is close enough to it so (\ref{eq:t34}) and (\ref{eq:t35}) hold.

Let $\{\mu_n\}$ be defined as in (\ref{10.30.07.2}) and $\{p_n\}$ be a sequence of monic polynomials such that $\nu(p_n) = \mu_n$. Further, let $\{q_n\}$ be a sequence of polynomials defined as at the beginning of the theorem when $\theta\in(0,1)$ and take $q_n$ to be an arbitrary polynomial of degree $n-k_n$ with zeros on $S_1$ when $\theta=1$.
Define
\[
(T^\delta_n h)(z) := p_n(z)q_n(z)\int_{L_\delta}\frac{h(\tau)}{p_n(\tau)q_n(\tau)}d\omega_\delta(\tau,z), \quad z\in E, \quad h\in\pl_n,
\]
Then $T_n^\delta$ is an operator from $\pl_n\cap H^\infty(U_\delta)$ to $C(E)$ such that
\begin{equation}
\label{eq:t311}
\|T_n^\delta h\|_E \leq \left(\frac{\|p_nq_n\|_E}{\min_{L_\delta}|p_nq_n|}\right)  \|h\|_{L_\delta}.
\end{equation}
Recall that
\[
\|p_nq_n\|_E^{1/n} \to \exp\left\{-\min_E U^{\mu_\theta+\lambda_\theta}\right\}=\exp\{-\widehat m_\theta\} \quad \mbox{as} \quad n\to\infty
\]
by (\ref{10.30.07.1}) and (\ref{10.30.07.2}), (\ref{10.25.07.1}), and since $E$ is regular. Moreover, the counting measures of zeros of $p_nq_n$, namely $\mu_n+\lambda_n$, are supported on $E\cup\Gamma$ and converge weakly to $\mu+\lambda$ that is supported on $E\cup S_\theta$. Therefore, we always can modify $q_n$, if needed, in such a manner that no zeros of $q_n$ lie in some neighborhood of $L_\delta$ and $\lambda_n$ still have the same asymptotic behavior. Hence, since the supports of $\mu_n+\lambda_n$ stay away from $L_\delta$, it holds that
\[
|p_nq_n|^{1/n} \to \exp\left\{-U^{\mu_\theta+\lambda_\theta}\right\} \quad \mbox{as} \quad n\to\infty \quad \mbox{uniformly on } \quad L_\delta.
\]
Thus, we get for the operator norm of $T_n^\delta$ that
\begin{eqnarray}
\limsup_{n\to\infty} \|T_n^\delta\|^{1/n} &\leq& \exp\left\{\max_{L_\delta}U^{\mu_\theta+\lambda_\theta}-\widehat m_\theta\right\} \nonumber \\ \label{eq:t31}
{} &=&\exp\left\{\max_{L_\delta}\left(U_D^{\lambda_\theta}(z)-g(z,\infty)\right)\right\} \leq \exp\left\{m_\theta+\delta\right\}
\end{eqnarray}
by (\ref{eq:300}) and (\ref{eq:t35}).

On the other hand, it holds that
\[
T_n^\delta h=h, \quad h(z_{j,n})=0, \quad h\in\pl_n,
\]
where $z_{1,n},\ldots,z_{k_n,n}$ are the zeros of $p_n$. Indeed, this holds because the ratio $h/p_nq_n$ is analytic in $U_\delta$ (including at infinity since $\deg(p_nq_n)=n$) and continuous on $L_\delta$. Let $\phi_1,\ldots,\phi_{k_n}$ be polynomials of degree at most $n$ such that $\phi_j(z_{i,n})=\delta_{ij}$, where $\delta_{ij}$ is the usual Kronecker symbol. Then for any $h\in\pl_n$ we have
\begin{equation}
\label{eq:t32}
(T_n^\delta h)(z) = h(z) - \sum_{j=1}^{k_n}h(z_j)\left(\phi_j(z)-(T_n^\delta\phi_j)(z)\right).
\end{equation}
Clearly, the sum on the right-hand sum of (\ref{eq:t32}) belongs to a $k_n$-dimensional subspace of $C(E)$ spanned by $\phi_j-T_n^\delta\phi_j$, $j=1,\ldots,k_n$. Hence,
\begin{equation}
\label{eq:t33}
d_{k_n}(A_n^\infty;C(E)) \leq \exp\{n\delta\}\|T_n^\delta\|
\end{equation}
by (\ref{eq:t34}). Combining (\ref{eq:t33}) with (\ref{eq:t31}), we get
\[
\limsup_{n\to\infty}\left(\frac1n\log d_{k_n}(A_n^\infty;C(E))\right) \leq m_\theta
\]
since $\delta$ was arbitrary. Thus, (\ref{eq:L31}) follows from (\ref{eq:t38}) and the last limit.
\end{proof}

\bibliographystyle{plain}
\bibliography{PrSYa08}

\end{document}